\documentstyle[12pt]{amsart}
\input{amssymb.sty}
\input xy
\xyoption{all}
\theoremstyle{plain}
  \newtheorem{thm}{Theorem}[section]
  \newtheorem{lem}[thm]{Lemma}
  
  \newtheorem{prop}[thm]{Proposition}
\theoremstyle{definition}

  \newtheorem{notation}{Notation\!\!}

\theoremstyle{remark}
  \newtheorem{rem}[thm]{Remark}
\newcommand{\bM}{\bar{M}}
\newcommand{\EmbM}{E^-_{\bar{M}}}
\newcommand{\EpbM}{E^+_{\bar{M}}}

\newcommand{\Amp}{{\operatorname{Amp}}}

\newcommand{\Coh}{{\operatorname{Coh}}}
\newcommand{\Cok}{\operatorname{Cok}}

\newcommand{\ext}{\operatorname{ext}}
\newcommand{\Ext}{\operatorname{Ext}}

\newcommand{\Hilb}{\operatorname{Hilb}}
\newcommand{\Hom}{\operatorname{Hom}}

\newcommand{\IIm}{\operatorname{Im}}
\newcommand{\Ker}{\operatorname{Ker}}

\newcommand{\NS}{\operatorname{NS}}

\newcommand{\Pic}{\operatorname{Pic}}
\newcommand{\PP}{{\bf{P}}}

\newcommand{\rk}{{\operatorname{rk}}}

\newcommand{\Spec}{\operatorname{Spec}}

\newcommand{\ZZ}{{\mathbb{Z}}}

\begin{document}
\title[Desingularization of some moduli of sheaves]
{Desingularization of some moduli scheme of
stable sheaves on a surface}
\author{Kimiko Yamada}
\email{kyamada@@math.kyoto-u.ac.jp}
\address{Department of mathematics, Kyoto University, Japan}
\maketitle
\begin{abstract}
Let $X$ be a nonsingular projective surface over an
algebraically closed field with characteristic $0$, and
$H_-$ and $H_+$ ample line bundles on $X$
separated by only one wall of type $(c_1,c_2)$.
Suppose the moduli scheme $M(H_-)$ of rank-two
$H_-$-stable sheaves with Chern classes $(c_1,c_2)$ is non-singular.
We shall construct a desingularization of $M(H_+)$ by using $M(H_-)$.
As an application, we consider whether
singularities of $M(H_+)$ are terminal or not
when $X$ is ruled or elliptic.
\end{abstract}
\section{Introduction}
Let $X$ be a projective non-singular surface over 
an algebraically closed field with characteristic $0$.
$H$ an ample line bundle on $X$. Denote by $M(H)$ 
the coarse moduli scheme of rank-two $H$-stable
sheaves with fixed Chern class
$(c_1,c_2)\in\NS(X)\times\ZZ$.
In this paper we think about singularity and
desingularization of $M(H)$ from the view of wall-crossing problem
of $H$ and $M(H)$.\par
Let $H_-$ and $H_+$ be ample line bundles on $X$
separated by only one wall of type $(c_1,c_2)$.
For a parameter $a\in(0,1)$, one can define the $a$-stability
of sheaves on $X$ and have the coarse moduli scheme $M(a)$ of
rank-two $a$-stable sheaves with Chern classes $(c_1,c_2)$.
Let $a_-$ and $a_+\in(0,1)$ be minichambers (see Section
\ref{section:desing} for details).
Assume $M_-=M(a_-)$ is non-singular;
one can find such $a_-$ when $X$ is ruled or elliptic for example.
In Section \ref{section:desing} we construct a desingularization
$\bar{\pi}_+: \bar{M}\rightarrow M_+$ of $M_+=M(a_+)$ by
using $M_-$ and wall-crossing methods.
In Section \ref{section:Kb-Kp} we 
calculate $K_{\bar{M}}- \bar{\pi}_+^* K_{M_+}$.
In Section \ref{section:eg} we apply it to consider whether
singularities of $M_+$ are terminal or not
when $X$ is ruled or elliptic.
%
%
%
%
\par
Here we mention related topics. 
About singularities of moduli spaces, Vakil \cite{Vak:Murphy}
shows that every singularity type of finite type over $\ZZ$
appears on moduli scheme of torsion-free sheaves on $\PP^5$,
and asks how about moduli scheme of sheaves on surfaces.
Thereby one can regard $M(H)$ as a model in which various
kinds of singularities can appear.
However a little is known about specific way to study
singularities of $M(H)$.
Methods in this article are suited to study 
what kind of singularities moduli scheme of sheaves has.
Perhaps one can use them to find interesting examples of singularity.
Properties of singularities in Section \ref{section:eg} seems to
relate with theory of determinantal variety 
over curve (see Remark \ref{rem:determ}).
This topic shall be studied in another article.
\begin{notation}
For a $k$-scheme $S$, $X_S$ is $X\times S$ and
$\Coh(X_S)$ is the set of coherent sheaves on $X_S$.
For $s\in S$ and $E_S\in\Coh(X_S)$, $E_s$ means $E\otimes k(s)$.
For $E$ and $F\in\Coh(X)$, 
$\ext ^i(E,F):=\dim\Ext^i_X(E,F)$ and $\hom(E,F)=\dim\Hom_X(E,F)$.
$\Ext_X^i(E,E)^0$ indicates
$\Ker(\operatorname{tr}:\Ext^i(E,E)\rightarrow H^i({\cal O}_X))$.
For $\eta\in\NS(X)$, we define $W^{\eta}\subset \Amp(X)$ by
$\bigl\{ H\in\Amp(X) \bigm| H\cdot\eta=0 \bigr\}$.
\end{notation}
\section{Desingularization of $M_+$ by using $M_-$ }\label{section:desing}
We begin with background material.
Let $H_-$ and $H_+$ be ample divisors lying in neighboring chambers
of type $(c_1, c_2)\in\NS(X)\times\ZZ$,
and $H_0$ an ample divisor in the wall $W$ of type $(c_1,c_2)$
which lies in the closure of chambers containing $H_-$ and $H_+$
respectively.
(Refer to \cite{Qi:equivalence} about the definition of wall and
chamber.) 
Assume that $M=H_+ -H_-$ is effective.
For a number $a\in [0,1]$ one can define the $a$-stability of a
torsion-free sheaf $E$ using
\[P_a(E(n))=\left. \{ (1-a)\chi(E(H_-)(nH_0))+ a\chi(E(H_+(nH_0)))
\}\right/ \rk(E).\]
There is the coarse moduli scheme ${\cal M}(a)$ of rank-two $a$-semistable
sheaves on $X$ with Chern classes $(c_1, c_2)$.
Denote by $M(a)$ its open subscheme of $a$-stable sheaves.
When one replace $H_{\pm}$ by $NH_{\pm}$ if necessary, 
$M(0)$ (resp. $M(1)$) equals
the moduli scheme of $H_-$-semistable (resp. $H_+$-semistable) sheaves.
There exist finite numbers $a_1\dots a_l\in(0,1)$ called minichamber
such that ${\cal M}(a)$ and $M(a)$ changes only when $a$ passes a
miniwall.
Refer to \cite[Section 3]{EG:variation} about these facts.
Fix numbers $a_-$ and $a_+$ separated by the only one miniwall,
and indicate ${\cal M}_{\pm}={\cal M}(a_{\pm})$ and
$M_{\pm}=M(a_{\pm})$ for short.
From \cite[Section 2]{Yam:Dthesis}, the subset
\begin{align*}
 {\cal M}_- \supset P_-= & \left\{ [E] \bigm| \text{$E$ is not 
$a_+$-semistable}\right\}  \\
 \bigl( \,\text{resp.}\quad {\cal M}_+ \supset P_+= & 
 \left\{ [E] \bigm| \text{$E$ is not $a_-$-semistable}\right\} \bigr) 
\end{align*}
is contained in $M_-$ (resp. $M_+$) and 
endowed with a natural closed subscheme structure of $M_-$
(resp. $M_+$).
Let $\eta$ be a element of
\[ A^+(W)=\left\{ \eta\in\NS(X) \bigm| \text{$\eta$ defines $W$,
$4c_2-c_1^2+\eta^2\geq 0$ and $\eta\cdot H_+>0$})\right\}. \]
After \cite[Definition 4.2]{EG:variation} we define
\[ T_\eta= M(1, (c_1+\eta)/2, n) \times M(1,(c_1-\eta)/2,m ),\]
where $n$ and $m$ are numbers defined by
\[ n+m=c_2-(c_1 -\eta^2)/4 \quad\text{and}\quad
   n-m= \eta\cdot(c_1-K_X)/2+(2a-1)\eta\cdot M,\]
and $M(1,(c_1+\eta)/2)$ is the moduli scheme of rank-one
torsion-free sheaves on $X$ with Chern classes $((c_1+\eta)/2,n)$.
If $F_{T_\eta}$ (resp. $G_{T_\eta}$) is the pull-back of a universal
sheaf of $M(1,(c_1+\eta)/2,n)$ (resp. $M(1,(c_1-\eta)/2,m)$) to
$X_{T_{\eta}}$, then we have an isomorphism
\[ P_- \simeq \underset{\eta\in A^+(W)}{\coprod} 
 \PP_{T_{\eta}}\left( Ext^1_{X_{T_{\eta}/T_{\eta}}}(F_{T_{\eta}},
G_{T_{\eta}}(K_X))\right) \]
from \cite[Section 5]{Yam:Dthesis}.
\begin{prop}[\cite{Yam:Dthesis} Proposition 4.9]\label{prop:Dthesis}
The blowing-up of $M_-$ along $P_-$ agrees with
the blowing-up of $M_+$ along $P_+$. So we have blowing-ups
\[ M_- \overset{\pi_-}{\longleftarrow} \tilde{M}
 \overset{\pi_+}{\longrightarrow} M_+.\]
\end{prop}
Now we assume that $P_-$ is nowhere dense in $M_-$ and that
every $E\in M_-$ satisfies that $\Ext_X^2(E,E)^0 =0$,
and explain how to induce a desingularization $\tilde{M}\rightarrow M_+$
of $M_+$ from $M_-$.
The idea is as follows. The problem of comparing $M_-$ and
$M_+$ is called wall-crossing problem or polarization-change
problem.
According to \cite{Yam:Dthesis}, we can endow a natural subset
\[ M_-\supset P=\bigl\{ [E] \bigm| \text{$E$ is not $H_+$-semistable}
\bigr\} \]
with a natural closed subscheme structure, and have a morphism
from the blowing-up of $M_-$ along $P$ to $M_+$
\[ M_- \longleftarrow \tilde{M}=B_P(M_-)
 \longrightarrow M_+ \]
in such a way that one can compare a universal family of $M_-$ with
that of $M_+$, if exists.
Since $M_-$ is non-singular, $\tilde{M}$ would be a desingularization
of $M_+$ if the center $P$ is non-singular.
We shall find a natural sequence of blowing-ups such that the
strict transform of $P$ becomes smooth, considering a relative
property of $P$ over $\Pic\times\Hilb\times\Hilb$.
As a result we obtain a diagram
\begin{equation*}
\xymatrix{
 & \bar{M} \ar[dl]_{\bar{\pi}_-} \ar[d] \ar[dr] \ar[dr]^{\bar{\pi}_+}
 & \\
 M_- & \tilde{M} \ar[l] \ar[r] & M_+,}
\end{equation*}
where $\bar{\pi}_-$ is a sequence of blowing-ups with smooth
centers, and so $\bar{M}$ is non-singular.\par
Now suppose that $A^+(W)=\{\eta\}$ for simplicity and denote $T_{\eta}= T$.
For a closed subscheme $Z\subset V$, let $B_Z(V)$ mean the blowing-up
of $V$ along $Z$.
First, in case where $\ext^1_X(F_t,G_t(K_X))$ is constant for all $t\in T$,
$P_-$ is non-singular and $B_{P_-}(M_-)$ is non-singular.
Hence from Proposition \ref{prop:Dthesis} we have a birational morphism 
$\pi_+:B_{P_-}(M_-)\rightarrow M_+$, which we can regard as a
desingularization of $M_+$.\par
In general case, set 
\begin{equation}\label{eq:l-1}
l_0=\min\{\ext^1(F_t,G_t(K_X))\bigm| t\in T\} \text{ and }
l_1=\max\{\ext^1(F_t,G_t(K_X))\bigm| t\in T\}.\end{equation}
Since one can readily show
$\ext^2(F_t,G_t(K_X))=\hom(G_t,F_t)=0$ for all $t\in T$,
there is a open covering $i: U\rightarrow T$ and a morphism
$F:V^0_U \rightarrow V^1_U$ of locally free ${\cal O}_U$-modules such
that $\rk V^1_U=l_1$, $\rk V^1_U-\rk V^0_U=-\chi(F_t,G_t(K_X))$, and
$P_- \sideset{}{_T}\times U \simeq {\bf P}_U(\Cok F)$.
Denote by $T_{l_1}\subset T$ the reduced closed subscheme with
\[T(k)\supset T_{l_1}(k)=
\{t\in T(k)\bigm| \dim\Ext^1(F_t,G_t(K_X))=l_1\}.\]
Let $C\subset T_{l_1}$ be any nonsingular subscheme  $C\subset T_{l_1}$
and denote $P_- \sideset{}{_T}\times C$ by $P_C$.
\begin{lem}\label{lem:outsideTl1}
We have an open covering $U' \rightarrow B_C(T)\setminus B_C(T_{l_1})$
and a homomorphism 
$F': V^0_{U'} \rightarrow V^1_{U'}$ such that
 $\IIm F'\supset \IIm F$, 
$\rk(\Cok F'\otimes k(t))\leq l_1 -1$ for all $t\in U'$, and
\[\PP_{U'}(\Cok F') \supset B_{P_C}(P_-) \sideset{}{_{B_C(T)}}\times U'. \]
\end{lem}
\begin{pf}
Assume $U=T$ for simplicity. 
Let $M(n,A)$ denote the square matrix ring of degree $n$ with coefficients in 
a ring $A$.
and $M(n,m,A)$ the set of matrixes of degree $(n,m)$.
Since $V^i_T$ is locally free, $F$ locally induces a matrix of degree
$(\rk V^0_T, \rk V^1_T)$.
From the definition of $T_{l_1}$, ideal sheaf $I_{T_{l_1}}\subset {\cal O}_T$
equals the radical of $(D_i)_i$, where $D_i$ runs over the set of all
minor determinant of degree $\rk V^1_T -l_1+1$.
Hence $B_C(T)\setminus B_C(T_{l_1})$ is covered by
$D_+(D_i)\subset B_C(T)=\operatorname{Proj}_T(\oplus I_C^d)$.
Let $A$ be a square matrix of degree $\rk V^1_T -l_1+1$ with
\begin{equation*}
F= \begin{pmatrix}
   A & * \\
   ** & *
   \end{pmatrix}.
\end{equation*}
Suppose $\det A\neq 0$ in ${\cal O}_T$ and $t\in T$ is the image
of a point $s\in D_+(\det A)$.
In the set of matrixes with coefficients in ${\cal O}_{T,t}$,
$A$ can be transformed into
\begin{math}
 \bigl(
 \begin{smallmatrix}
  I & 0 \\
  0  & A'
 \end{smallmatrix}
\bigr)
\end{math}
with a square matrix $A'$ with coefficients in $m_{T,t}$.\par
Since $C$ and $T$ are non-singular, we have a splitting
$T\simeq C\times S$ after replacing $T$ by an \'{e}tale neighborhood of
$t$. 
The point $t\in T$ corresponds to $(t_1,t_2)\in C\times S$.
Let $p_1:C\times S\rightarrow C$ and $p_2:C\times S\rightarrow S$
denote projections and
$i_1:C \simeq C\times \{t_2\} \hookrightarrow C\times S$ 
and $i_2:S \simeq \{t_1\}\times S \hookrightarrow C\times S$ 
immersions.
For $a\in m_{T,t}$, we have $i_{1\,*} p_{2\, *} i_{2\, *} a=0$
and so $p_{2\, *} i_{2\, *} a\in \Ker i_{1\, *}=I_{C,t}$.
Hence if $\tilde{A}'$ is the cofactor matrix of $A'$, then
the matrix $B= p_{2\, *} i_{2\, *}(\tilde{A}')/ \det A'$ is an element
of $M(d,{\cal O}_{B_C(T),s})$ because of the assumption on $A'$ and
the definition of $B_C(T)$.
Elements of the matrix $A'B-I$ lie in $\Ker i_{2\,*}=I_{T,t}$, and so
$A'$ is invertible in $M(d,{\cal O}_{B_C(T),s})$.
As a result 
$A$ is invertible in the $M(\rk V^1_T-l_1+1,{\cal O}_{B_C(T),s})$
and we obtain the lemma.
%
\end{pf}
\begin{lem}\label{lem:smooth}
The natural immersion 
$B_{P_C}(P_-) \sideset{}{_{B_C(T)}}\times B_C(T_{l_1})\subset 
P_-\sideset{}{_T}\times B_C(T_{l_1})$ is isomorphic, and
the projection 
$B_{P_C}(P_-) \sideset{}{_{B_C(T)}}\times B_C(T_{l_1}) \rightarrow B_C(T_{l_1})$
is smooth; that is a $\PP^{l_1-1}$-bundle.
\end{lem}
\begin{pf}
When we denote by $D\subset B_C(T_{l_1})$ the exceptional divisor, we have
\begin{equation*}
\xymatrix{
P_{T_{l_1}} \setminus P_C 
\simeq B_{P_C}(P_-) \sideset{}{_{B_C(T)}}\times [B_C(T_{l_1})\setminus D]
 \ar@{=}[r] \ar@{^{(}->}[d] &
 P_- \sideset{}{_T}\times [B_C(T_{l_1})\setminus D] 
\ar@{^{(}->}[d] \\
B_{P_C}(P_-) \sideset{}{_{B_C(T)}}\times B_C(T_{l_1}) 
\ar@{^{(}->}[r] 
& P_- \sideset{}{_T}\times B_C(T_{l_1}). }
\end{equation*}
 The first row is isomorphic, and $B_C(T_{l_1})\setminus D$
 is dense in $B_C(T_{l_1})$.  
 The second row is set-theoretically bijective since
 the dimension of fiber of
 $B_C(T_{l_1}) \sideset{}{_T}\times P_-\rightarrow B_C(T)$
 is upper-semicontinuous.
  The left-hand side of second row is reduced since $T_{l_1}$ is reduced
 so the second row is isomorphic.
\end{pf}
%
There is a sequence of blowing-ups
$T_{l_1}^{(N)}\rightarrow \dots \rightarrow T^{(1)}_{l_1} 
   \rightarrow T_{l_1}=T_{l_1}^{(0)} $
with nonsingular center $C^{(i)}\subset T_{l_1}^{(i)}$ such that
$T_{l_1}^{(N)}$ is non-singular.
They induces a sequence of blowing-ups
$T^{(N)}\rightarrow \dots \rightarrow T^{(1)}\rightarrow  T=T^{(0)}$.
%
Put $P^{(0)}=P_-$ and $P^{(1)}=B_{P_{C^{(0)}}}(P_-)$.
Then we have a natural morphism 
$P^{(1)}\rightarrow T^{(1)}=B_{C^{(0)}}(T)$ such that
the restriction
$P^{(1)} \sideset{}{_{T^{(1)}}}\times C^{(1)}\rightarrow C^{(1)}$
is smooth from Lemma \ref{lem:smooth}.
%
Set $P^{(2)}=B_{P^{(1)} \sideset{}{_{T^{(1)}}}\times C^{(1)}} (P^{(1)})$.
It is naturally a scheme over $T^{(2)}=B_{C^{(1)}}(T^{(1)})$, and
$P^{(2)} \sideset{}{_{T^{(2)}}}\times C^{(2)} \rightarrow C^{(2)}$ is
smooth from Lemma \ref{lem:smooth}.
Denote $M^{(0)}=M_-$, $M^{(1)}=B_{P_{C^{(0)}}}(M^{(0)})$ and
$M^{(2)}= B_{P^{(1)} \sideset{}{_{T^{(1)}}}\times C^{(1)}} (M^{(1)})$.
$P^{(i)}$ are closed subschemes of $M^{(i)}$ for $i\leq 2$.
%
In the same way we obtain a sequence of blowing-ups
$M^{(N)}\rightarrow\dots\rightarrow M^{(1)}\rightarrow M_-=M^{(0)}$
with centers $P^{(i)}\subset M^{(i)}$, which is smooth over $C^{(i)}$
for $0\leq i\leq N-1$.
The subscheme $P^{(N)} \sideset{}{_T}\times T_{l_1}^{(N)}$ of $P^{(N)}$
is smooth over $T_{l_1}^{(N)}$ by Lemma \ref{lem:smooth},
$T_{l_1}^{(N)}$ is non-singular by the definition,
and hence
$P^{(N)} \sideset{}{_{T^{(N)}}}\times T_{l_1}^{(N)}$ itself is smooth.
%
Let $M^{(N+1)}$ denote
$B_{P^{(N)} \sideset{}{_{T^{(N)}}}\times T_{l_1}^{(N)} } (M^{(N)})$.
Its closed subscheme
$P^{(N+1)}:= B_{P^{(N)} \sideset{}{_{T^{(N)}}}\times T_{l_1}^{(N)} }(P^{(N)})$
is smooth over a nonsingular scheme
$T^{(N+1)}:=B_{T_{l_1}^{(N)}}(T^{(N)})$ from Lemma \ref{lem:smooth}.
Using Lemma \ref{lem:outsideTl1},
we can find a homomorphism 
$F_0:V^0_{T^{(N+1)}} \rightarrow V^1_{T^{(N+1)}}$ such that
$\rk \Cok F_0 \otimes k(t)\leq l_1 -1$ for all $t\in T^{(N+1)}$
and that the $T^{(N+1)}$-scheme
$P^{(N+1)}$ is contained in $\PP_{T^{(N+1)}}(\Cok F_0)$
if we replace $T^{(N+1)}$ with an open covering.
%
Then repeat this process after changing
$T$ to $T^{(N+1)}$, $P_-$ to $P^{(N+1)}$, $M_-$ to 
$M^{(N+1)}:= B_{P^{(N)} \sideset{}{_{T^{(N)}}}\times 
T_{l_1}^{(N)}}(M^{(N)})$,
$F$ to $F_0$ and $l_1$ to $l_1 -1$.
Consequently we obtain a sequence of blowing-ups 
$T^{(N')}\rightarrow \dots \rightarrow T^{(0)}=T$ 
with non-singular center $C^{(i)}\subset T^{(i)}$
and blowing-ups
$M^{(N')}\rightarrow\dots \rightarrow M^{(0)}=M_-$ with
center $P^{(i)}\subset M^{(i)}$ as follows.
There is a commutative diagram
\begin{equation*}
\xymatrix{ P^{(i)} \ar[r] \ar[d] & T^{(i)}\ar[d] \\
P_- \ar[r] & T}
\end{equation*}
and a homomorphism $F^{(m)}:V^0_{T^{(N')}} \rightarrow V^1_{T{(N')}}$
as follows. It holds that $P_- \subset \PP_{T^{(N')}}(\Cok F^{(N')})$,
$\rk\Cok F^{(N')}\otimes k(t)\leq l_0$ for all $t\in t^{(N')}$,
and the first row of
\begin{equation*}
\xymatrix{ P^{(N')} \ar@{^{(}->}[r] \ar[d] & \PP_{T^{(N')}}(\Cok F^{(N')}) \ar[d]\\
P_- \ar@{=}[r] & \PP_T (\Cok F)}
\end{equation*}
is isomorphic when it is restricted to
the inverse image of
$T^{(N')}\setminus T^{(N')}\sideset{}{_T}\times T_{l_0+1}$.
Thus $P^{(N')}\rightarrow T^{(N')}$ is smooth by the same proof as
Lemma \ref{lem:smooth}.
Set $M^{(N'+1)}=B_{P^{(N')}}(M^{(N')})$. Since $M^{(N'+1)}\rightarrow M_-$ is
a composition of blowing-ups of the smooth scheme $M_-$
along nonsingular centers, $M^{(N'+1)}$ itself is nonsingular.
One can verify that $M^{(N'+1)}\rightarrow M_-$ splits as
$M^{(N'+1)}\rightarrow B_{P_-}(M_-)$ , so
we obtain a morphism 
$M^{(N'+1)}\rightarrow B_{P_-}(M_-)=B_{P_+}(M_+)\rightarrow M_+$ 
from Proposition \ref{prop:Dthesis},
which is a desingularization of $M_+$ since $M^{(N'+1)}$ is
nonsingular.
Put $\bar{M}=M^{(N'+1)}$. We have constructed
\begin{equation}\label{eq:desing}
\xymatrix{
  & \bar{M} \ar[dl]_{\bar{\pi}_-} \ar[dr]^{\bar{\pi}_+} 
  \ar[d]_{\tilde{\pi}} &  \\
 M_-  & \tilde{M} \ar[l]^{\pi_-} \ar[r]_{\pi_+} & M_+} 
\end{equation}
with $\bar{M}$ nonsingular.
\section{Calculation of  $K_{\bar{M}}- \bar{\pi}_+^*K_{M_+}$ }\label{section:Kb-Kp}
Assume that 
$M_+\supset \operatorname{Sing}(M_+):=\{E \bigm| \ext^2(E,E)^0\neq 0 \}$
satisfies $\operatorname{codim}(M_+, \operatorname{Sing}(M_+))\geq 2$
and that $P_+\subset M_+$ is nowhere dense, and hence both $M_-$ and
$M_+$ are locally complete intersections and so Gorenstein schemes.
%
Let us calculate $K_{\bar{M}}-\bar{\pi}_+^* K_{M_+}$.
By construction, each step $M^{(i+1)}\rightarrow M^{(i)}$ in
$\bar{M}\rightarrow M_-$ is a $\PP^{l^{(i)}-1}$-bundle over
$C^{(i)}\subset T^{(i)} \sideset{}{_T}\times T_{l^{(i)}} \subset T^{(i)}$,
where $l_0\leq l^{(i)}\leq l_1$ and $T_l\subset T$ is the reduced subscheme
such that 
$T_{l^{(i)}}=\{t\in T \bigm| \dim\Ext^1 (G_t, F_t)\geq l^{(i)}\}$
set-theoretically.
If we denote by $D^{(i)}\subset \bar{M}$ the pull-back of
the exceptional divisor of $M^{(i)}\rightarrow M^{(i-1)}$, then
\begin{equation}\label{eq:Kb-Km}
K_{\bar{M}} - \bar{\pi}_-^* K_{M_-}=
\sum_{i=0}^{N'-1} \operatorname{codim}(P^{(i)}, M^{(i)})\, D^{(i)}=
\sum_i [\dim M_- -(l^{(i)}-1+\dim C^{(i)})]\, D^{(i)}.
\end{equation}
Next consider $\bar{\pi}_-^*(K_{M_-})- \bar{\pi}_+^*(K_{M_+})$.
By the proof of Proposition \ref{prop:Dthesis}, which uses elementary
transform, we have the following.
\begin{prop}\label{prop:Dthesis-2}
Denote the exceptional divisor $\pi_-^{-1}(P_-)=\pi_+^{-1}(P_+)$
by $D$.
Suppose we have a universal family $E^-_{M_-}\in\Coh(X_{M_-})$
of $M_-$ and a universal family $E^+_{M_+}\in\Coh(X_{M_+})$ of $M_+$.
If $\pi:D\rightarrow P_+\rightarrow T$ is a natural map,
then there are $T$-flat modules $F_T$ and $G_T$ on $X_T$,
line bundles $L_{\pm}$ on $P_{\pm}$ and
a line bundle $L_0$ on $\tilde{M}$ such that we have exact sequences
\begin{equation}\label{eq:elemtrans}
 0 \longrightarrow \pi_+^* E^+_{M_+} \otimes L_0 \longrightarrow
 \pi_-^* E^-_{M_-} \overset{f}{\longrightarrow} 
 \pi^*G_T \otimes \pi_+^*L_+  \longrightarrow 0
\end{equation}
in $\Coh(X_{\tilde{M}})$ and
\begin{equation}\label{eq:elemtrD}
 0 \longrightarrow \pi^*F_T\otimes \pi_-^*L_- \longrightarrow
 \pi_-^*(E^-_{M_-})|_{X_D} \overset{f|D}{\longrightarrow}
 \pi^*G_T\otimes \pi_+^*L_+ \longrightarrow 0
\end{equation}
in $\Coh(X_D)$.\\
\end{prop}
The exact sequence \eqref{eq:elemtrD} is the relative $a_+$-Harder
Narashimhan filtration of $E^-_{M_-}$.
Here we remark that
generally a universal family of $M_-$ exists \`{e}tale-locally, but
one can generalize this proposition to general case with straightforward labor.
Suppose $L_{\pm}$ and $L_0$ in this proposition are trivial for
simplicity.
From \eqref{eq:elemtrans}
\begin{align*}
\bar{\pi}_+^* K_{M_+} & - \bar{\pi}_-^* K_{M_-} \\
 = & \pi_-^*\det {\bf R}Hom_{X_{M_-}/M_-}(E^-_{M_-},E^-_{M_-})-
\pi_+^*\det {\bf R}Hom_{X_{M_+}/M_+}(E^+_{M_+},E^+_{M_+}) \\
 =&
\det {\bf R}Hom_{X_{\bM}/\bM}(\pi_-^*E^-_{M_-}, \pi_-^*E^-_{M_-})-
\det {\bf R}Hom_{X_{\bM}/\bM}(\pi_+^*E^+_{M_+}, \pi_+^*E^+_{M_+}) \\
 =&
\det {\bf R}Hom_{X_{\bM}/\bM}(\EmbM, \EpbM)+
\det {\bf R}Hom_{X_{\bM}/\bM}(\EpbM, \pi^* G_T) \\
 & -\det {\bf R}Hom_{X_{\bM}/\bM}(\EmbM, \EpbM)+
\det {\bf R}Hom_{X_{\bM}/\bM}(\pi_+^* G_T, \EpbM)\\
 = & 
\det {\bf R}Hom_{X_{\bM}/\bM}(\EmbM, G_D)+ 
\det {\bf R}Hom_{X_{\bM}/\bM}(G_D, \EpbM).
\end{align*}
If $i:D\hookrightarrow M_-$ is inclusion, then by \eqref{eq:elemtrD}
\begin{multline}\label{eq:Hom-EmG}
\det{\bf R}Hom_{X_{\bM}/\bM}(\EmbM, G_D)=
\det i_*{\bf R}Hom_{X_D/D}(\EmbM |_D, G_D)= \\
\det i_*{\bf R}Hom_{X_D/D}(F_D,G_D)+ \det i_*{\bf R}Hom_{X_D/D}(G_D,G_D).
\end{multline}
Since $\det {\cal O}_D=D$, \eqref{eq:Hom-EmG} equals
$ [\chi(F_t,G_t)+\chi(G_t,G_t)]\,D$ for any $t\in D$.
By the Serre duality
\begin{multline*}
\det{\bf R}Hom_{X_{\bM}/\bM}(G_D, \EpbM)=
\det{\bf R}Hom_{\bM}( {\bf R}Hom_{X_{\bM}/\bM}(\EpbM, G_D(K_X)),
{\cal O}_{\bM})\\
= -\det {\bf R}Hom_{X_{\bM}/\bM}(\EpbM, G_D(K_X))=
-\det i_*{\bf R}Hom_{X_D/D}(\EpbM |_D, G_D(K_X)) \\
= -[\chi(F_t,G_t(K_X))+\chi(G_t,G_t(K_X))]\,D=
-[\chi(G_t,F_t)+\chi(G_t,G_t)]\,D.
\end{multline*}
Therefore
$\bar{\pi}_+^* K_{M_+} - \bar{\pi}_-^* K_{M_-}=[\chi(F_t,G_t)-\chi(G_t,F_t)]D$.
On the other hand, we put
\begin{equation}\label{eq:lambda}
\tilde{\pi}^* D=\sum_{i=0}^{N'} \lambda_i\,D^{(i)}.
\end{equation} 
If $\lambda_i$ is determined, then
we can calculate $K_{\bar{M}}-\bar{\pi}_+^* K_{M_+}$ by
\eqref{eq:Kb-Km} and \eqref{eq:lambda}.
Let $Z_{M_+}\subset P_+$ denote the pull-back of 
$\sideset{}{_{i=l_0+1}^{l_1}}\cup \operatorname{Sing}(T_i)\subset T$ by 
$P_+\rightarrow T$, which is a nowhere-dense closed subscheme.
Let us consider the induced open subset 
\begin{equation}\label{eq:Up}
U_{M_+}=M_+\setminus Z_{M_+}.
\end{equation}
One can regard $M_-\supset P_-\supset P^{(0)}$ \'{e}tale-locally as
$ k[x_1,\dots,x_m]\supset I_{P^{(0)}}=(x_1,\dots,x_n)\supset 
I_{P_-}=(f_1,\dots,f_{n'})$.
We have 
\[ B_{P^{(0)}}(M_-)_{(x_n)}=\Spec k[x_1/x_n=x'_1,\dots,x_{n-1}/x_n=x'_n,x_n,
 \dots,x_m],\]
and if 
$\pi_0:M^{(1)}:=B_{P^{(0)}}(M_-)\rightarrow M_-$ is a natural morphism, then
\[k[x'_1,\dots,x'_{n-1},x_n,\dots,x_m]\supset {\pi}_0^{-1}I_{P_-}\cdot 
{\cal O}_{P^{(1)}}= 
\bigl( x_n^{N_1}\bar{f}_1(x'_i,x_j),\dots,
x_n^{N_{n'}}\bar{f}_{n'}\bigr)\]
where $\bar{f}_{i}$ cannot be divided by $x_n$.
$\lambda_0$ equals $\max (N_i)_i$.
It becomes $1$ when $\dim M_->l_1-1+\dim C^{(0)}$.
Indeed for any point $t\in P^{(0)}$,
$\rk \Omega_{P^(0)}\otimes k(t)\leq l_1-1+\dim C^{(0)}$ since
$P^{(0)}$ is a $\PP^{l_1-1}$-bundle over $C^{(0)}$.
Hence in the exact sequence 
\[ CN_{P_-/M_-}\otimes k(t) \overset{\tau}{\longrightarrow}
  \Omega_{M_-}\otimes k(t) \longrightarrow 
  \Omega_{P_-}\otimes k(t) \longrightarrow 0, \]
$\tau$ cannot be zero if $\dim M_->l_1-1+\dim C^{(0)}$,
so $N_j=1$ for some $j$.
Over $U_{M_+}$ we can suppose $T_j$ is non-singular, and
so $P^{(i)} \sideset{}{_{M_+}}\times U_{M_+}$ does not contain
any irreducible component of the exceptional divisor of
$M^{(i)}\rightarrow \bar{M}$
when $\dim T_l< \dim T_{l-1}$ for all $l_0<l\leq l_1$.
Thereby, similarly to the case where $i=0$, one can show that
$\lambda_i\geq 1$, and that if $\dim M_-> l^{(i)}-1+\dim C^{(i)}$
then $\lambda_i=1$, since
$P^{(i)}$ is a $\PP^{l^{(i)}}$-bundle over $C^{(i)}$
for $l_0\leq l^{(i)}\leq l_1$.
Thus we have shown the following.
\begin{prop}\label{prop:Kb-Kp}
In the diagram \eqref{eq:desing} it holds that
\begin{multline}\label{eq:Kb-Kp}
 K_{\bar{M}}- \bar{\pi}_+^* K_{M_+} \\
 \sum_{i=0}^{N'} [\dim M_- -(l^{(i)}-1+\dim C^{(i)})+\lambda_i
\{ \chi(F_t,G_t)-\chi(G_t,F_t)\}]\, D^{(i)}.
\end{multline}
Suppose $\dim T_l<\dim T_{l-1}$ for all $l_0<l\leq l_1$ and
$\dim M_-> l_1-1+\dim T$. If the image of $T^{(i)}\subset M^{(i)}$
in in $T$ agrees with $T_j$ for
some $j$, then $\lambda_i=1$.
\end{prop}
Remark that the image of $D^{(i)}$
in $T$ agrees with $T_j$ for some $j$ if
$D^{(i)}$ has non-empty intersection with
$\pi_+^{-1}(U_{M_+})$.
Thus one can use this proposition to verify whether singularities
in $U_{M_+}$ are terminal or not.
\begin{rem}\label{rem:determ}
When the image of $D^{(i)}$ in $T$ does not agree with $T_j$
for any $j$, the value $\lambda_i$ seems to relate with
determinantal varieties over $C$.
\end{rem}
\section{Examples: ruled or elliptic surface}\label{section:eg}
We shall give examples of $M_+$ with $M_-$ non-singular.
If a morphism $X\rightarrow C$ to a nonsingular curve $C$ exists, then
by \cite[p.142]{Frd:holvb} we have a $(c_1,c_2)$-suitable polarization,
that is, an ample line bundle $H$ such that $H$ does not lie on
any wall of type $(c_1,c_2)$, and for any wall $W=W^{\eta}$ of
type $(c_1,c_2)$, we have $\eta\cdot f=0$ or
$\operatorname{Sign}(f\cdot \eta)=\operatorname{Sign}(H\cdot \eta)$.
From \cite[p.159, p.201]{Frd:holvb}, if $X$ is a ruled surface
or an elliptic surface, then
any rank-two sheaf $E$ of type $(c_1,c_2)$ which is stable respect to
$(c_1,c_2)$-suitable polarization is good, i.e.
$\Ext^2(E,E)^0=0$.\par
%
%
(A) First we suppose that $X$ is a (minimal) ruled surface.
When $c_1\cdot f$ is odd $M(H)$ is empty for $(c_1,c_2)$-suitable
polarization. Thus we assume $c_1=0$.
If a rank-two sheaf $E$ of type $(c_1,c_2)$ is stable with respect
to a polarization $H$ such that $H\cdot K_X<0$, then $E$ is good
and so $M(H)$ is nonsingular.
Hence we assume that $W^{K_X}\cap \Amp(X)\neq\emptyset$, so
$2\leq g=g(C)$ and $e(X)\leq 2g-2$ from the
description of $\Amp(X)$ \cite[Prop. V.2.21]{Ha:text}.
Since $\dim \NS(X)=2$,
if we move polarization $H$ from a $(c_1,c_2)$-suitable one,
then $M(H)$ may begin to admit singularities when $H$ passes the
wall $W^{K_X}$.
Let $H_-$ and $H_+$ be ample line bundles separated by only one wall
$W^{K_X}$.
$M(H_-)$ is non-singular, and $E^+\in \PP_+$ has a non-trivial
exact sequence
\begin{equation}\label{eq:NSF-wkx}
0 \longrightarrow G=L\otimes I_{Z_l} \longrightarrow E^+
  \longrightarrow F=L^{-1}\otimes I_{Z_r} \longrightarrow 0
\end{equation}
with $-2L\sim mK_X$.
About this filtration we have $\Ext^2_-(E^+,E^+)=0$ since $p_g(X)=0$
(See \cite[p. 49]{HL:text} for $\Ext_{\pm}$), and
\[\ext^2(E^+,E^+)=\ext^2_+(E^+,E^+)=
 \ext^2(L\otimes I_{Z_l},L^{-1}\otimes I_{Z_r})=
 hom(I_{Z_r},{\cal O}(K_X+2L)I_{Z_l}).\]
Since $W^{K_X}$ defines a wall, $H^0({\cal O}(K_X+2L))=0$
unless $2L+K_X=0$.
Hence $\ext^2(E^+,E^+)^0\neq 0$ if and only if $-2L=K_X$ and 
$Z_l \subset Z_r$.
As a result
when one defines $a$-stability using $H_{\pm}$, 
\[\chi^a(\otimes I_{Z_r}, -K_X\otimes I_{Z_l})=Aa+B+l(Z_l)\]
for some constant $A$ and $B$, and so
the moduli scheme $M(a)$ of $a$-stable sheaves begins to
admit singularities just when $a$ passes a miniwall $a_0$
defined by
\begin{equation*}
 l(Z_l)=\left\{
 \begin{gathered}
 c_2/2-(g-1)\quad\text{if $c_2$ is even} \\
 (c_2-1)/2-(g-1) \quad\text{if $c_2$ is odd.}
 \end{gathered}\right.
\end{equation*}
Let $a_-$ and $a_+$ be minichambers separated by only one
miniwall $a_0$.
$M(a_+)=M_+$ has singularities along $P_+ \sideset{}{_T}\times T'$,
where 
\[T'=\bigl\{ (L\otimes I_{Z_l}, L^{-1}\otimes I_{Z_r}) \bigm| -2L=K_X
 \bigr\}_{red} \subset
M(1,K_X/2,l(Z_l))\times M(1,-K_X/2,l(Z_r)).\]

%
%
(B) Suppose that $X$ is an elliptic surface with a section $\sigma$
and $c_1=\sigma$.
In contrast to ruled surfaces, $K_X^2=0$ and so $W^{K_X}\cap\Amp(X)$ is
always empty, though one can study some singularity appearing in $M(H)$
by Proposition \ref{prop:Kb-Kp}.
Let $\pi: X\rightarrow C$ be an elliptic fibration,
$f\in \NS(X)$ its fiber class,
$d=- \deg R^1\pi_*({\cal O}_X) -\sigma^2 \geq 0$.
We have a natural map to a ruled surface 
$\kappa:X\rightarrow \PP(\pi_*({\cal O}(2\sigma)))=\PP({\cal E}_2)$.
Since $\kappa_*(\sigma)$ is a section of 
$\PP({\cal E}_2)$, and since the pull-back of an ample line bundle by
a finite map is ample,
$L=a f$ satisfies $W^{2L-c_1}\cap \Amp(X)\neq\emptyset$ if 
$a>0$ from the description of the ample cone of a ruled surface.
Let $c_1$ be $\sigma$ and
$c_2=(c_1-L)\cdot L=a$.
Then any sheaf $E$ with non-trivial exact sequence
\begin{equation}\label{eq:Ext-ellipt}
 0 \longrightarrow F=L \longrightarrow E \longrightarrow G=L^{-1}\otimes c_1
\longrightarrow 0,
\end{equation}
whose Chern class equals $(c_1,c_2)$,
is stable with respect to a $(c_1,c_2)$-suitable ample line bundle.
Indeed, $(2L-c_1)\cdot f<0$ and so
$\pi_*({\cal O}(2L-c_1))=0$ and $R^1\pi_*({\cal O}(2L-c_1))$
commutes with base change.
Thus the exact sequence
\[ 0 \longrightarrow H^1(C,\pi_*({\cal O}(2L-c_1)))) \longrightarrow 
  H^1(X,{\cal O}(2L-c_1))
 \longrightarrow H^0(E,R^1\pi_*({\cal O}(2L-c_1))) \]
shows that the restriction of the exact sequence \eqref{eq:Ext-ellipt}
to a general fiber is non-trivial, and so
a corollary of
Artin's theorem for vector bundles on an elliptic curve 
\cite[p. 89]{Frd:holvb} and a basic property of a suitable polarization
\cite[p. 144]{Frd:holvb} deduce that $E$ is stable with respect to
a suitable polarization. Thereby such $E$ is good.
Let $H_-$ and $H_+$ be ample line bundles which lie in no wall
of type $(c_1,c_2)$
with $(2L-c_1)\cdot H_-<0< (2L-c_1)\cdot H_+$.
One can define $a$-stability by them.
Let $a_0$ be a miniwall such that 
$\chi^{a_0}({\cal O}(L))=\chi^{a_0}({\cal O}(2L-c_1))$,
$a_-<a_0<a_+$ minichambers, and $M_{\pm}=M(a_{\pm})$.
Then some connected components of $P_-\subset M_-$ contains
any sheaf $E$ with non-trivial exact sequence \eqref{eq:Ext-ellipt},
and some neighborhood of them in $M_-$ is non-singular.
It induces a desingularization of some open neighborhood of
connected components ${\cal K}_+$
of $P_+$ consisting of sheaves $E^+$
with a non-trivial exact sequence
\[ 0 \longrightarrow L^{-1}\otimes c_1 \longrightarrow E^+ \longrightarrow
   L \longrightarrow 0 \]
as in Section \ref{section:desing}.\par
We have in case of (A) $\ext^1(G,F)\leq 1$, and in case of (B) 
$\ext^1(G,F)=h^0(c_1-2L+K_X)-\chi(c_1-2L)\leq 2c_2+C(X)$
with some constant $C(X)$ independent of $c_2$ because
$h^0(c_1-2L+K_X)=0$ if $a=c_2$ is sufficiently large.
Thus in both cases one can show that, if $c_2$ is sufficiently
large, then all singularities of $M_+$ along 
the dense open set $U_{M_+}\cap P_+$  in 
$P_+\subset \operatorname{Sing}(M_+)$ defined at
\eqref{eq:Up} are terminal.
\providecommand{\bysame}{\leavevmode\hbox to3em{\hrulefill}\thinspace}
\providecommand{\MR}{\relax\ifhmode\unskip\space\fi MR }
 \MRhref is called by the amsart/book/proc definition of \MR.
\providecommand{\MRhref}[2]{%
  \href{http://www.ams.org/mathscinet-getitem?mr=#1}{#2}
}
\providecommand{\href}[2]{#2}

\end{document}